\documentclass{amsart}
\usepackage{amsmath,amssymb,enumerate}

\usepackage{epsfig,fancyhdr,color}

\newtheorem{theorem}{\rm\bf Theorem}
\newtheorem{proposition}{\rm\bf Proposition}[section]
\newtheorem{lemma}[proposition]{\rm\bf Lemma}
\newtheorem{corollary}[proposition]{\rm\bf Corollary}

\theoremstyle{definition}
\newtheorem{definition}[proposition]{\rm\bf Definition}

\theoremstyle{remark}
\newtheorem{remark}[proposition]{\rm\bf Remark}

\newtheorem{example}[proposition]{\rm\bf Example}

\def\interieur#1{\mathord{\mathop{\kern 0pt #1}\limits^\circ}}

\title{Weak metrics on Euclidean domains}

\author{Athanase Papadopoulos}
\address{A. Papadopoulos, Institut de Recherche Math{\'e}matique Avanc\'ee,
Universit{\'e} Louis Pasteur and CNRS,
7 rue Ren\'e Descartes,
 67084 Strasbourg Cedex - France} \email{papadopoulos@math.u-strasbg.fr}

\author{Marc Troyanov}
\address{M. Troyanov, Section de Math{\'e}matiques,  \'Ecole Polytechnique F{\'e}derale de
Lausanne, 1015 Lausanne - Switzerland}
\email{marc.troyanov@epfl.ch}

\date{February 10, 2006}

\begin{document}

\begin{abstract}
A \emph{weak metric} on a set is a function  that satisfies the axioms of a metric except the symmetry and the separation axioms. In the present paper 
 we introduced a weak metric, called the  \emph{Apollonian weak metric}, on any subset of a Euclidean space which is either bounded or whose boundary is unbounded. We then relate  this weak metric  to some familiar metrics such as  the  Poincar\'e metric, the Klein-Hilbert metric, Funk  metric, and the part metric 
which  play  an important role in classic and  recent work on geometric function theory.

\bigskip



\end{abstract}

\maketitle

\section{introduction}\label{s1}

The axioms for a metric space were formulated exactly 100 years ago, in a famous paper by Maurice Fr\'{e}chet, see \cite{Frechet}. 
Since then, several important  generalizations of the notion of metric space
appeared.  In the present paper, we shall consider weak metrics and semi-metrics. We first recall the definitions.

\begin{definition}[Weak metric and semi-metric]   A {\it weak metric} on a set $X$ is a function
$\delta:X\times X\to [0,\infty)$   satisfying
\begin{enumerate}[i)]
 \item $\delta(x,x)=0$ for all $x$ in $X$;
\item  $\delta(x,y)+\delta(y,z)\geq \delta(x,z)$ for
all $x$, $y$ and $z$ in $X$.
\end{enumerate}

A {\it semi-metric} is a symmetric weak metric, that is, a weak metric satisfying 
\begin{enumerate}[i)] \setcounter{enumi}{2}
 \item  $\delta(x,y)=\delta(y,x)$ for all $x$ and $y$
in $X$.
\end{enumerate}
\end{definition}

H. Busemann studied extensively functions 
satisfying some of the axioms of a metric (see \cite{Busemann1944, Busemann1970}), and he called them ``general metrics" or simply
 ``metrics". The name ``weak metric" is due to H. Ribeiro (\cite{Rib}).
There are two distinct notions of separation for weak metric, namely:
\begin{definition}
A weak metric $\delta$ is said to be {\it strongly separating} if we have 
$$\min{\{\delta(x,y),\delta(y,x)\}} = 0 \quad \Leftrightarrow  \quad  x=y,$$
and it is said to be {\it weakly separating} if
$$\max{\{\delta(x,y),\delta(y,x)\}} = 0 \quad \Leftrightarrow  \quad  x=y,$$
for all $x$ and $y$ in $X$.
\end{definition}

\medskip

In the case of a semi-metric, these two notions clearly coincide. A \emph{metric} on $X$, in the usual sense, is a separating semi-metric. 

\medskip

Given a  weak metric $\delta:X\times X\to [0,\infty)$, it  can be interesting  to consider a
{\it symmetrization} of it, and to try to compare this symmetrization with other known metrics or semi-metrics.  In fact, there exist several notions  of
symmetrization, none of them being more natural than the others. We
shall deal in this paper with two of these, defined as follows.

\begin{definition}[Symmetrizations] \label{symm} Let $\delta:X\times X\to[0,\infty)$ be a 
weak metric. A {\it symmetrization} 
of $\delta$ is  one of the following functions
 $\sigma \delta:X\times X\to [0,\infty)$ and  $S \delta:X\times X\to [0,\infty)$,
defined for $x$ and $y$ in $X$ by  
$$\sigma \delta(x,y)=\max\left\{\delta(x,y),\delta(y,x)\right\}$$
and 
$$S \delta(x,y) = \frac{1}{2}\left(\delta(x,y)+\delta(y,x)\right).$$
Both symmetrizations are semi-metrics.
The semi-metric  $\sigma \delta$ is sometimes called the \emph{max-symmetrization} and  $S \delta$ the \emph{meanvalue-symmetrization} of $\delta$.
\end{definition}

Observe the following inequalities:
$$
 S \delta(x,y) \leq \sigma \delta(x,y) \leq 2 S\delta(x,y).
$$
Furthermore, the equality $S \delta = \sigma \delta$ holds  if and only if $\delta$ is itself symmetric. In that case, both symmetrizations coincide with $\delta$.

It is also clear from the definitions  that both  symmetrizations of a weakly separating weak metric gives rise to a metric. 

\medskip

In order to be more concrete, we start right away with an example.
Let $A$ be an open subset of $\mathbb{E}^n$, with $A\not=\mathbb{E}^n$.
We introduce the function ${i}_A$ on $A\times A$ defined by
$$
{i}_A(x,y) =  \log \left( 1 + \frac{|x-y|}{d(x,\partial A)}\right)
= \sup_{a\in \partial A}   \log \left( 1 + \frac{|x-y|}{|x-a|}\right)
$$
for all $x,y\in A$.

\medskip

\begin{proposition} The function $i_A$ is a weak metric on $A$.
\end{proposition}

\proof  We prove the triangle inequality. For $x$, $y$ and $z$ in $A$, we have
$$
\vert y-z\vert\geq \vert x-z\vert -\vert x-y\vert
$$
and 
$$\vert x-y\vert + d(x,\partial A)\geq d(y,\partial A).$$
Multiplying the two inequalities, we obtain
$$
\vert y-z\vert   \left(\vert x-y\vert + d(x,\partial A)\right)  \geq \left(\vert x-z\vert -\vert x-y\vert\right) d(y,\partial A)
$$
or, equivalently,
$$
\left(d(x,\partial A)+\vert x-y\vert\right)\left(d(y,\partial A)+\vert y-z\vert\right)\geq
\left(d(x,\partial A)+\vert x-z\vert\right)d(y,\partial A).
$$

The last inequality is equivalent to
$$\left(1+{ \vert x-y\vert\over d(x,\partial A)}\right)
\left(1+ {\vert y-z\vert\over d(y,\partial A)}\right)\geq 
\left(1+ {\vert x-z\vert\over d(x,\partial A)}\right).$$ Taking logarithms, we obtain the triangle inequality for 
$i_A$.\qed

\medskip

In their study of uniform domains in Euclidean spaces \cite{GO}, 
Gehring and Osgood considered a metric which is the symmetrization $S\i_A$ of our weak metric $i_A$. Its is usually denoted
by $\tilde{j}_A$ and it is therefore defined by
$$
\tilde{j}_A(x,y) =  S \mathrm{i}_A(x,y)  =
\frac 1{2} \left\{ \log \left( 1 + \frac{|x-y|}{d(x,\partial D)}\right)
+\log  \left( 1 + \frac{|x-y|}{d(x,\partial D)}\right)  \right\}.
$$

\bigskip

The   symmetrization  $\sigma i_A$ of this weak metric is a metric that 
has been considered by M. Vuorinen in his study of 
conformal invariants and moduli of families of curves, see \cite{Vuorinen}.
Vuorinen's metric is usually written as  
$$
j_A(x,y)  =
 \log \left( 1 + \frac{|x-y|}{\min \{\delta(x), \delta(y)\}}\right)
$$
where $\delta(z) = d(z,\partial D)$.

\medskip

There is a large literature on the metrics $\tilde{j}_A$ and $j_A$, see for instance \cite{Seittenranta} and \cite{Hasto2006}.

\medskip  

As the above example illustrates,  we shall see in this paper that a certain number of important metrics are naturally obtained as a symmetrization of some weak metric. The
weak metric appears then as a kind of primitive structure on which the actual metric is built. It is then an interesting question to investigate the
 geometric properties of the weak metric and to compare
them with those of the associated symmetric metric (or semi-metric).
To our knowledge, this question  has not been really studied so far. Let us formulate it as the following general:

\medskip

\textbf{Problems.}  (1) \emph{Given a metric space $(X,d)$, find a natural weak metric $\delta$ on $X$ such that $d$ is the symmetrization of $\delta$, i.e. $d= S\delta$ or $d= \sigma \delta$.} \\
(2) \emph{Describe the geometry of $(X,\delta)$ and compare it to the geometry of $(X,d)$.}

\medskip

By the word geometry, we mean here the study of geodesics, of isometries, of curvature and so on.
\medskip

These problems are not precisely formulated, in particular we should not expect question (1) to have a unique answer. 
However, we believe that these questions are worth 
investigating, at least 
in the case of some important metric spaces.

\medskip
In this paper, we shall in particular address these problems in the case of the hyperbolic plane. We shall see that the hyperbolic metric can be obtained as a symmetrization of at least three natural weak metrics. 
The first weak metric is the so called \emph{Funk weak metric} 
and it is related to the projective (Klein) model of hyperbolic geometry. The other two are the \emph{Apollonian weak metrics},
 and they are related to the conformal (Poincar\'{e}) model of the unit disk and  the upper-half plane
 respectively. We shall give explit formulas for the Apollonian weak metrics of the upper-half plane and the unit disk in Theorem \ref{Th1} and Theorem
  \ref{variation-S1} respectively. We 
 shall observe that the isometry group of the Apollonian weak metric is quite different from the isometry group of the hyperbolic metric. On the other hand, we 
 shall show in Theorem \ref{Th3} that the hyperbolic lines in the unit disk are geodesics for the Apollonian weak metric.

\medskip

The  weak metrics that we consider here are defined for a wide class of domains in  Euclidean space.
The Funk weak metric is classical although not so popular. The Apollonian weak metric is a new
notion which we define in section \ref{s3} below. The name ``Apollonian" was chosen because a
symmetrization of this weak metric is the Apollonian semi-metric discussed in the paper \cite{Beardon1995} by Beardon.

\medskip

The plan of the rest of this paper is as follows.
In Section \ref{s2}, we recall the definition of the Funk weak metric, and we identify its two symmetrizations. In Section \ref{s3}, we define  the Apollonian weak metric. and we relate it to some other known metrics. 
In Section \ref{s4}, we give an explicit  formula for 
the Apollonian weak metric of the unit disk in $\mathbb{C}$ and we draw some consequences of that formula. The Poincar\'{e} metric of the unit disk is a symmetrization of the Apollonian weak metric in the same way as the Klein-Hilbert metric of that disk is a symmetrization of the Funk metric.
Section \ref{s5} is concerned about  the notion of geodesic associated to a weak metric. We study in particular the geodesics of the Apollonian weak metric of the unit disk in $\mathbb{C}$.

\medskip
We shall denote by $\vert 
x-y\vert$ the Euclidean distance between the points $x$ and $y$ in a Euclidean space.

\section{The Funk weak metric}\label{s2}

The \emph{Funk weak metric}  is a weak metric defined on open bounded convex subsets of $\mathbb{E}^n$. It was discovered by P. Funk \cite{Funk} 
and is discussed in  \cite{Busemann1970,Zaustinsky}. To recall its definition, 
let $A$ be a nonempty bounded open convex domain in $\mathbb{E}^n$. 
For every $x$ in $A$, we set  $\mathcal{F}(x,x)=0$. For $y$ distinct from $x$, 
we consider the  Euclidean ray starting at $x$ and passing through $y$, we let 
 $a$ be the intersection point
of  that ray with the boundary of $A$, and we set
\begin{equation}
 \mathcal{F}(x,y)=\log{\frac{|x-a|}{|y-a|}}.
\end{equation} 
The function  $\mathcal{F}$ is a strongly separating weak metric on $A$, which we call the Funk weak metric.  The proof of the triangle inequality follows from the theorem 
 of Menelaus (see \cite[ Appendix I]{Zaustinsky}).

\medskip

If we need to indicate the dependence of the Funk weak metric with respect to the domain, we denote it by $\mathcal{F}_A$. 

\medskip

We recall that a {\it similarity} defined on a domain in $\mathbb{E}^n$  is a map $\phi$ satisfying 
$\vert \phi(x)-\phi(y)\vert=\mu\vert x-y\vert$ for some $\mu>0$ and for all $x,y$ in ths domain.  A similarity always extends
 as a global affine map $\phi : \mathbb{E}^n \to \mathbb{E}^n$. Similarities form a subgroup of the affine group.
Observe that  
any similarity $\phi : A \to B$ between bounded convex open subsets $A,B \subset \mathbb{E}^n$  is an isometry from $(A,\mathcal{F}_A)$
 to $(B,\mathcal{F}_B)$. This is obvious
since similarities preserve ratios of Euclidean distances.

\medskip

The mean-value symmetrization (in the sense of Definition \ref{symm}) of the Funk weak metric give rise to the \emph{Klein-Hilbert metric}. 
Let us  describe this metric  $\mathcal{H}$. \label{Hilbert}.  Let $A$ be again a 
nonempty open bounded convex subset of $\mathbb{E}^n$.  
For $x=y\in A$, we set $\mathcal{H}(x,y)=0$, and for $x\not=y$, we consider the Euclidean line containing $x$ and $y$. It intersects the boundary of $A$ in two points, $a$ and $b$, these names chosen such that 
$b,x,y,a$ follow each other in that order on that line. We then set
\begin{equation}\label{hilbert}
\mathcal{H}(x,y)= \frac{1}{2} \log\left( {\vert x-a\vert\over \vert y-a\vert}  {\vert y-b\vert\over \vert 
x-b\vert}\right).
\end{equation}
Note that the quantity $\displaystyle {\vert x-a\vert\over \vert y-a\vert}  {\vert y-b\vert\over \vert  x-b\vert}$ is equal to the cross ratio $[b,x,y,a]$ of the four points.
Since the notions of Euclidean line and of cross ratio of aligned points are invariant under projective transformations, the Klein-Hilbert metric is also
 invariant under projective transformations  (which include the similarities).

This metric was first defined by F. Klein on the unit disk $\mathbb{D}^2\subset\mathbb{E}^2$. 
It defines what is usually called the Klein 
model of hyperbolic geometry on the disk $\mathbb{D}^2$. This metric was later on defined by D. Hilbert, using the same formula, on an 
arbitrary bounded open covex subset of  $\mathbb{E}^n$.
Formula (\ref{hilbert}) shows that the Klein-Hilbert metric is a
symmetrization of the Funk semi-metric. More precisely, we have
\begin{equation}\label{symF}
 \mathcal{H}= S\mathcal{F}.
\end{equation}

\section{The part metric} 

This notion was introduced by H. S. Bear in his study of complex function
algebras, \cite{Bear1991,BB}. It can be defined in the following abstract setting. Consider a set $X$ and a class  $\mathcal{B}$ of positive
real-valued functions on $X$. We
 introduce an equivalence relation (which we call ``part-equivalence") on $X$ as follows: two points $x$ and $y$ in $X$ are equivalent if and only if there exists a constant $c>0$ such that the following 
 Harnack-type inequality
$$
 \frac{1}{c} \leq \frac{u(y)}{u(x)}  \leq c
$$
holds for all functions in $\mathcal{B}$. The equivalence classes are called the corresponding \emph{parts} of $X$. 
They form a partition of $X$ which is associated to $\mathcal{B}$. Such partitions were
 first considered by A. M. Gleason.

\medskip

On every part of $X$, we have the following natural metric:
\begin{equation}\label{p1}
 p(x,y) = p_{\mathcal{B}}(x,y) = \sup \left\{\left\vert \log \left({\frac{u(x)}{u(y)}} \right) \right\vert \; \big| \, u \in \mathcal{B} \right\},
\end{equation}
which is called the \emph{part metric} induced from $(X,\mathcal{B})$.

\medskip

\begin{proposition}
 Let $A$ be an open bounded convex subset of $\mathbb{E}^n$ and let 
 $\mathcal{B}$ be the class of positive  functions on $A$ that are 
 restrictions of affine functions $u : \mathbb{E}^n \to \mathbb{R}$. 
 Then, all the elements of $A$ are part-equivalent for the relation
 induced by $\mathcal{B}$ the corresponding part metric on $A$ 
 coincides with the max-symmetrization of the Funk weak metric:
\begin{equation*} \label{ }
p_{\mathcal{B}}(x,y) =  \sigma \mathcal{F} (x,y) 
 =\max\left\{\log{\vert x-a\vert\over\vert y-a\vert},
\log{\vert y-b\vert\over\vert x-b\vert}\right\}.
\end{equation*}
\end{proposition}

\proof  That all the points of $A$ are part-equivalent will follow from the fact that 
$p(x,y)$ is finite for all $x$ and $y$ in $A$, which follows from what we prove
now.  Given $x,y \in A$ we denote by $a$ and $b$  the two points  lying on the intersection of 
the boundary $\partial A$ and the Euclidean line passing through $x$ and $y$, assuming $b,x,y,a$ 
follow each other in that order on the line.
To prove that $\mathcal{F} (x,y) \leq  p_{\mathcal{B}}(x,y)$, we consider
an affine function $u$ such that $u > 0$ on $A$
and $u(a) = 0$. For $t>0$, set $z(t) = ty + (1-t)a$. Then $y = z(1) $ and $x = z(s)$ for  
$s = \frac{|x-a|}{|y-a|}$. Furthermore, there exists $\lambda>0$ such that
for $t>0$ we have $u(z(t)) = \lambda t$, since $u(z(0)) = u(a) = 0$ and $u(x) > 0$. 
Thus
$$
 p(x,y) \geq \log \frac{u(x)}{u(y)} = \log \frac{\lambda s}{\lambda} =
 \log{\frac{|x-a|}{|y-a|}} = \mathcal{F}(x,y).
$$
A similar argument shows that $p(x,y) \geq  \mathcal{F}(y,x)$ and thus 
$p \geq \sigma  \mathcal{F}$.

To prove the converse inequality, we consider an arbitrary  affine function $v$ such that $v > 0$ on $A$. We parametrize
the segment $[a,b]$ by  $z(t) = tb + (1-t)a$ \ so that 
$v(z(t)) = \lambda t + \mu$ for some $\lambda, \mu$.
For $0<t\leq s$, we have
$$
 \frac{v(z(t))}{v(z(s))} = \frac{\lambda t + \mu}{\lambda s + \mu}
 =  \frac{t + \mu/\lambda }{s + \mu/\lambda }
 \leq \max \left\{\frac{s}{t} , \frac{t}{s} \right\}.
$$
It easily follows from this inequality that $p(x,y) \leq 
\max \left\{\mathcal{F}(y,x), \mathcal{F}(x,y) \right\}.$

\qed

\bigskip

Let us now consider an open subset $A$ of the complex plane $\mathbb{C}$ on which there
exists a non-constant positive harmonic function (for this it suffices that $\hbox{Card}(\mathbb{C}\setminus A)\geq 2$), and let  $\mathcal{B}$
be the set of harmonic functions on $A$. 
The corresponding part metric is thus given by 
\begin{equation}\label{p.harm}
 p(x,y)  = \sup \left\{\left\vert \log \left({\frac{u(x)}{u(y)}} \right) \right\vert \; \big| \,   u > 0, \text{harmonic in } A \right\},
\end{equation}
Since the composition of a positive harmonic function with a conformal transformation is again a positive harmonic function, the part metric is a conformal invariant of domains. 

\medskip

It is worthwile to note that in the case where $A$ is the unit disk $\mathbb{D}^2$,   the 
corresponding metric space  is essentially isometric to the hyperbolic plane. More precisely, H. Bear proved in \cite[Corollary 1]{Bear1991}
the following
\begin{proposition}
In the unit disk $\mathbb{D}^2$, the part metric (\ref{p.harm}) associated to the class of harmonic functions coincides with twice the Poincar\'e  metric of that disk.
\end{proposition}

\qed

\section{The Apollonian weak metric}\label{s3}

 Let $A\subset\mathbb{E}^n$ be an open subset
 and let  $\partial A=\overline{A}\setminus A$ be its boundary. In this section, we suppose that either
 $A$ is bounded or  $\partial A$ is unbounded. Note that any nonempty convex subset $A$ of $\mathbb{E}^n$ with $A\not=\mathbb{E}^n$ satisfies
 theis hypothesis.

 We define a function $\delta_A$ on $A\times A$ by the formula
\begin{equation}\label{delta-a}
\delta_A(x,y)=\sup_{a\in\partial A}\log {\vert x-a\vert\over\vert y-a\vert}.
\end{equation}

\begin{proposition}\label{unbounded} The function $\delta_A$ is a weak metric. 
\end{proposition}

\proof The proof is straightforward. We only say a few words on the fact that  $\delta_A$ is nonnegative.
First, suppose that
$A$ is bounded. For any distinct points $x$ and $y$ in $A$, 
consider the Euclidean ray starting at $x$ and
passing through $y$ and 
let $a$ be an intersection point of that ray with $\partial A$.
We
have $\displaystyle {\vert x-a\vert \over \vert y-a\vert}>1$, which implies $\delta_A(x,y)\geq 0$. 
Now suppose that $\partial A$ is unbounded. Let $(a_n)$ be a 
sequence of points in $\partial A$  such that $\vert x-a_n\vert \to \infty$ for some (or equivalently, for any) $x$ in $A$.
 Then, for any $x$ and $y$ in $A$, we have $\displaystyle {\vert x-a_n\vert \over \vert y-a_n\vert} \to 1$ as $n  \to\infty$, which shows $\delta_A(x,y) \geq 0$. 

\qed

 \medskip
 \begin{definition}[The Apollonian weak metric]\label{Apo}
 For any open subset  $A\subset\mathbb{E}^n$  which is either bounded  or whose boundary $\partial A$ is unbounded,
 the weak metric provided by Proposition \ref{unbounded} is called the {\it Apollonian weak} metric of $A$. (The name is chosen because of Definitions
 \ref{Beardon} and \ref{Hasto}, and
 Proposition \ref{sym} below.)
 \end{definition}

 \medskip

 The following invariance property is straightforward.

\begin{proposition}
For
 any similarity $\phi$ of $\mathbb{E}^n$, we have, 
for every $x$ and $y$ in $A$,
$$\delta_A(x,y)=\delta_{\phi(A)}(\phi(x),\phi(y)).$$
\end{proposition}

\qed

We have the following easy comparison between the Apollonian weak metric $\delta_A$ and the weak metric $i_A$ that we defined in \S\ref{s1}:

\begin{proposition}\label{dual}
For every $x$ and $y$ in $A$, we have $\delta_A(x,y)\leq i_A(y,x)$.
\end{proposition}

\begin{remark}
 Observe that in the statement above, the last term is $i_A(y,x)$ and not $i_A(y,x)$. We note that for any given weak metric $\delta$ the function
 $\delta'$ defined by $\delta'(y,x)=\delta(x,y)$ is also a weak metric, which can be called the weak metric \emph{dual} to $\delta$. Therefore, we can consider  Proposition 
 \ref{dual} as giving a comparison between the weak metric  $\delta_A$ and the weak 
 metric dual to $i_A$.
\end{remark}

\proof For  $x$ and $y$ in $A$ and for $z$ in $\partial A$, we can write
$$
\log {\vert x-z\vert \over \vert y-z\vert} \leq  \log {\vert x-y\vert +\vert y-z\vert \over \vert y-z\vert}
=\log \left(1+{\vert y-x\vert\over \vert y-z\vert}\right).
$$
Taking the supremum over $z\in\partial A$, we obtain the desired result.
\qed

\medskip

We shall see examples of weak metrics on sets $A$ satisfying both kinds of hypotheses of Definition \ref{Apo}. 
The first example is the following:

 \begin{example}[The Apollonian weak metric on the upper half-plane] \label{ex-halfplane}
In this example,  $A$ is the upper half-plane  $\mathbb{H}^2 = \left\{ \left.
z \in \mathbb{C} \hspace{0.25em} \right| \hspace{0.25em}
\text{Im} ( z ) > 0 \right\}$. The associated Apollonian weak metric $\delta_{\mathbb{H}^2}:\mathbb{H}^2\times \mathbb{H}^2\to [0,\infty)$ is given by
\begin{equation*}
  \delta_{\mathbb{H}^2}(x, y)=  \sup_{a \in \mathbb{R}} \log  \frac{|x-
a|}{|y - a|}  .
\end{equation*}
It is easy to see that the  restriction of $\delta_A$  to the 
vertical half-line   $\{z = is\}$, $s>0$ is given by
\begin{equation}\label{def.delta}
  \delta_{\mathbb{H}^2}(x,y)  = \max \left\{0,\log {t\over s} \right\}
   =
   \left\{ \begin{array}{cll}
     \left. \log {t\over s} \right.   & \text{if} &s \leq t
  \\
     0 & \text{if} &s \geq t
   \end{array} \right.
\end{equation}
for $x=is$ and $y=it$.
From this, we can see that the weak metric $\delta_A$  is  not symmetric and   not strongly separating. It is weakly separating.
This Apollonian weak metric has been studied  in the paper \cite{BPT}, in  with a weak metric introduced by Thurston on the 
Teichm\"uller space of a hyperbolic surface \cite{Thurston1985} and it was shown that $\delta_A$ coincides with Thurston's geometrically defined weak metric,
 if we interpret the upper half-plane as the Teichm\"uller space of Euclidean metrics on the torus.
The following result was obtained in \cite[Proposition 3]{BPT}:
\begin{theorem}\label{Th1}
The Apollonian weak metric of the upper-half plane is given by 
 \begin{equation}
 \delta_{\mathbb{H}^2}=\log\left({\vert y-\overline{x}\vert+\vert y- {x}\vert\over
 \vert x-\overline{x}\vert}\right)
 \end{equation}
 for every $x$ and $y$ in $\mathbb{H}^2$.
\end{theorem}

\qed

It follows that the symmetrization $S\delta_{\mathbb{H}^2}$ coincides with
the Poincar\'e metric $h_{\mathbb{H}^{2}}$ on  $\mathbb{H}^2$. In other words, for all
$x$ and $y$ in  $\mathbb{H}^2$, we have
\begin{equation}\label{poincmet}
 S\delta_{\mathbb{H}^2}(x,y) = h_{\mathbb{H}^{2}}(x,y)=\frac{1}{2} \log \left(
 \frac{|x - \bar{y} | + |x - y|}{|x -\bar{y}|
  - |x - y|} \right).
\end{equation}
\end{example}

Next, we want to relate  the Apollonian weak metric
to some semi-metrics that appear in recent works 
of Beardon and others.  We first recall these semi-metrics.

\begin{definition}[The Apollonian semi-metric]\label{Beardon}
Let
$A$ be any open subset of $\mathbb{E}^n$.
The {\it Apollonian semi metric} $\alpha_A$ on $A$ is defined by 
$$\alpha_A(x,y)=\sup_{a\in\partial A}\log{ {\vert x-a\vert\over \vert y-a\vert}} + 
\sup_{b\in\partial A}\log{ {\vert y-b\vert\over \vert
x-b\vert}}=\sup_{a,b\in\partial A}[b,x,y,a].$$
\end{definition}

The name ``Apollonian" was given by A. Beardon who studied that semi-metric in  \cite{Beardon1995} This semi-metric is also discussed in several later papers,
for example \cite{GH}.
(As mentionned in \cite{GH}, this metric has been earlier introduced by D. Barbilian \cite{barbilian},  and Beardon rediscovered it independently). The Apollonian semi-metric is a metric in the usual sense of the word
  if  $\partial A$ does not lie in an $(n-1)$-dimensional sphere 
or an $(n-1)$-dimensional hyperplane (see \cite[ Theorem 1.1]{Beardon1995}). The Apollonian semi-metric is  invariant under M\"obius transformations.  
If $A$ is  the upper-half plane or 
 the unit ball of $\mathbb{E}^n$, the  Apollonian semi-metric coincides with the Poincar\'e metric of these spaces.

\begin{definition}[The half-Apollonian semi-metric]\label{Hasto} Let
$A$ be an open subset of   $\mathbb{E}^n$.
The {\it half-Apollonian semi-metric} $\eta_A$ on $A$ is defined by 
$$\eta_A(x,y)=\sup_{a\in\partial A}\left\vert \log{ {\vert x-a\vert\over \vert 
y-a\vert}}\right\vert.$$
\end{definition}

The half-Apollonian semi-metric was introduced by P. H\"ast\"o,  \& H. Lind\'en in   
\cite{Hasto-Linden}. It is invariant under similarities (cf. \cite[Theorem 1.2]{Hasto-Linden} )
and it is a metric in the usual sense of the word whenever $\mathbb{E}\setminus A$ is
not contained in a hyperplane of $\mathbb{E}^n$.

\medskip

The following proposition, whose proof is immediate from the definitions, shows the relation between the Apollonian and the half-Apollonian semi-metrics and the Apollonian weak metric $\delta_A$:

\begin{proposition}\label{sym} for any open subset $A$ of $\mathbb{E}^n$,
 we have
$$\sigma\delta_A=\eta_A \qquad \text{and} \qquad S\delta_A=\alpha_A.$$
\end{proposition}

\qed

\section{The Apollonian weak metric of the unit disk}\label{s4}

The rest of the paper is mainly devoted to a discussion of the  Apollonian weak metric in the 
unit disk $ \mathbb{D}^{2}\subset \mathbb{C}$.  In this section, we give an explicit formula for that
weak metric.

\medskip

\begin{theorem}\label{variation-S1}
The Apollonian weak metric $\delta_{\mathbb{D}{}^{2}}$  is given by the following
formula:
\begin{equation}\label{formulaD2}
\delta_{\mathbb{D}{}^{2}}(x,y)=\log\left(\frac{\left|x-y\right|+\left|x\overline{y}-1\right|}{\left|1-\left|y\right|^{2}\right|}\right).
\end{equation}
\end{theorem}

\proof The result follows directly from the first statement of Proposition \ref{S-1} below.

\qed

\medskip

Before stating the needed Proposition, we first draw a few consequences of formula (\ref{formulaD2}).

\begin{corollary} \label{co-sy}
The symmetrization $S\delta_{\mathbb{D}^{2}}$ of the weak metric $\delta_{\mathbb{D}^{2}}$ 
on the unit disk  $\mathbb{D}^{2}$ coincides
with the Poincar\'{e} metric $h_{\mathbb{D}^{2}}$ of that disk:
\[
S\delta_{\mathbb{D}{}^{2}}(x,y) = {1\over 2}\left(\delta_{\mathbb{D}{}^{2}}(x,y)+\delta_{\mathbb{D}{}^{2}}(y,x) \right)=  h_{\mathbb{D}^{2}}=
{1\over 2}\log\left(\frac{\left|1-x\overline{y}\right|+\left|x-y\right|}{\left|1-x\overline{y}\right|-\left|x-y\right|}\right).
\]
\end{corollary}

\proof The proof is a direct calculation from Theorem \ref{variation-S1}. 
Observe that the result also follows from Proposition \ref{sym} and the
result of Beardon stating that the Apollonian semi-metric of the unit disk is the Poincar\'e
metric. 

\qed

\medskip

\begin{corollary}\label{co-non}
The Apollonian weak metric $\delta_{\mathbb{D}{}^{2}}$  
is nonsymmetric, unbounded and weakly separating. 
\end{corollary}

\proof  Using Formula (\ref{formulaD2}), we obtain the following special
values:
$$\delta_{\mathbb{D}{}^{2}}(x,0)=\log\left|1+|x|\right|, \qquad
\delta_{\mathbb{D}{}^{2}}(0,x) = -\log\left|1-|x|\right|$$
Thus  $\delta_{\mathbb{D}{}^{2}}$ is non-symmetric and unbounded since, $\delta_{\mathbb{D}{}^{2}}(0,x)\to\infty$
as $\vert x\vert\to 1$. The fact that it is weakly separating follows from Corollary \ref{co-sy}.

\medskip

\qed

\medskip

\begin{corollary}
The Apollonian weak metric $\delta_{\mathbb{D}{}^{2}}$  
is not invariant under the group of M\"obius transformation preserving the unit disk.
\end{corollary}

This property is in contrast with a property of the hyperbolic metric.

\proof   Given an arbitrary pair of points $x,y \in \mathbb{D}^2$, there exists a M\"obius transformation 
$g$ preserving the disk and exchanging $x$ and $y$ ($g$ is the  $180^0$-hyperbolic rotation around
the mid-point of the hyperbolic segment $[x,y]$). If $\delta_{\mathbb{D}{}^{2}}$ were invariant under
the M\"obius group, then we would have $\delta_{\mathbb{D}{}^{2}}(y,x) = \delta_{\mathbb{D}{}^{2}}(g(x),g(y)) = \delta_{\mathbb{D}{}^{2}}(x,y)$,
which contradicts Corollary \ref{co-non}.

\qed

\bigskip

The following Proposition was used in the proof of Theorem  \ref{variation-S1}. For later use, we formulate a more complete statement  than what is needed in that proof.

\begin{proposition}\label{S-1}
Let us fix two distinct points $x$ and $y$ in $\mathbb{C}$, and consider the  
function $f : \mathbb{S}^{1}\rightarrow [0,\infty)$ defined by 
\[
f(a)=\left|\frac{x-a}{y-a}\right|.
\]

The maximum value of this function on the circle $\mathbb{S}^{1}$ is given by 
\[
\max_{|a|=1}f(a)=\frac{\left|x-y\right|+\left|x\overline{y}-1\right|}{\left|\left|y\right|^{2}-1\right|}
\]
and this maximum is achieved at a unique point $a^{+}(x,y)\in \mathbb{S}^{1}$
given by
 \[
a^{+}(x,y)=\frac{\left|x-y\right|(x\overline{y}-1) y+(x-y)\left|x\overline{y}-1\right|}{\left|x-y\right|(x\overline{y}-1)+(x-y)\left|x\overline{y}-1\right|\overline{y}}.
\]

\smallskip

The minimum of $f$ on $\mathbb{S}^{1}$ is given by 
\[
\min_{|a|=1}f(a)=\left|\frac{\left|x-y\right|-\left|x\overline{y}-1\right|}{\left|y\right|^{2}-1}\right|
\]
and it is achieved at a unique point 
\[
a^{-}(x,y)=\frac{\left|x-y\right|(x\overline{y}-1) y-(x-y)\left|x\overline{y}-1\right|}{\left|x-y\right|(x\overline{y}-1)-(x-y)\left|x\overline{y}-1\right|\overline{y}}.
\]
\end{proposition}

\medskip

To prove the proposition, we shall use the following lemma:

\medskip

\begin{lemma}\label{lemma}
Let $g:\mathbb{C}\to\mathbb{C}$ be the function given by 
\[
g(z)=\lambda(\mu z+1),
\]
where $\lambda,\mu\in\mathbb{C}$ and $\mu\not=0$. 
Then $$\max_{|z|=1}|g|=|\lambda|(|\mu|+1)$$ and this maximum is achieved
at the unique point $z^{+}= |\mu|/\mu$. 
Likewise $$\min_{|z|=1}|g|=|\lambda|\left||\mu|-1\right|$$ and this
minimum is achieved at the unique point $z^{-}=-|\mu|/\mu$.
\end{lemma}

\medskip

\proof Without loss of generality, we can assume $\lambda=1$. Then, for $\vert z\vert=1$, the point $g(z)=\mu z+1$ in the complex plane describes, as $z$ varies,  a circle of centre 1 and radius $\vert
\mu\vert$. The function $\vert g(z)\vert$ is the distance from that point to the origin, and therefore it has a unique maximum, which is equal to $1+\vert
\mu\vert$. The
rest of the proof follows by an analogous reasoning.

\qed 

\bigskip

{\it Proof of Proposition \ref{S-1}.} Let us set
 $\displaystyle q= {\overline{y}a-1\over y-a}$. then, we have $\displaystyle  a = {qy+1\over q+\overline{y}}$. 
We claim that $\vert a\vert=1\iff \vert q\vert=1$. Indeed, if $\vert a\vert=1$, we can write $a=e^{i\theta}$ with $\theta\in\mathbb{R}$. Then, 
$\displaystyle 
q=\frac{\overline{y}e^{i\theta}-1}{y-e^{i\theta}}= - e^{i\theta}\frac{e^{-i\theta}-\overline{y}}{e^{i\theta}-y}$
which shows that $\vert q\vert=1$. In the same way, we can see that if $\vert q\vert=1$ then $\vert a\vert=1$.

Now set
\begin{eqnarray*}
g(q) = \frac{x-a}{y-a} & = &\left(x-\left(\frac{qy+1}{q+\overline{y}}\right)\right).\left(y-\left(\frac{qy+1}{q+\overline{y}}\right)\right)^{-1}\\
 & = & \left(\frac{x\overline{y}-1}{|y|^{2}-1}\right)\cdot \left(\left(\frac{x-y}{x\overline{y}-1}\right) q+1\right)\cdot
 \end{eqnarray*}

Applying Lemma \ref{lemma} with $\displaystyle \lambda= \frac{x\overline{y}-1}{|y|^{2}-1} $
and $\displaystyle \mu= \frac{x-y}{x\overline{y}-1}$, we see that 
\begin{eqnarray*}
\max_{|a|=1}f(a)  =  \max_{|q|=1}|g(q)| & = & \left|\frac{x\overline{y}-1}{|y|^{2}-1}\right|\cdot\left(\left|\frac{x-y}{x\overline{y}-1}\right|+1\right)\\
 & = & \frac{\left|x-y\right|+\left|x\overline{y}-1\right|}{\left|\left|y\right|^{2}-1\right|}.\end{eqnarray*}
 \\
This maximum is achieved at the unique point
\[
q^{+}=\left|\frac{x-y}{x\overline{y}-1}\right|\cdot\left(\frac{x-y}{x\overline{y}-1}\right)^{-1}=\frac{(x\overline{y}-1)|x-y|}{|x\overline{y}-1|(x-y)},
\]
which corresponds to
 \[
a^{+}=\frac{q^{+}y+1}{q^{+}+\overline{y}}=\frac{\left|x-y\right|(x\overline{y}-1) y+\left|x\overline{y}-1\right|(x-y)}
{\left|x-y\right|(x\overline{y}-1)+\left|x\overline{y}-1\right|(x-y)\overline{y}}.
\]
Likewise, we have
\begin{eqnarray*}
\min_{|a|=1}f(a)  =  \min_{|q|=1}|g(q)| & = & \left|\frac{x\overline{y}-1}{|y|^{2}-1}\right|\cdot\left|\left|\frac{x-y}{x\overline{y}-1}\right|-1\right|\\
 & = & \left|\frac{\left|x-y\right|-\left|x\overline{y}-1\right|}{\left|y\right|^{2}-1}\right|
 \end{eqnarray*}
and this minimum is achieved at the unique point
 \[
q^{+}=\left|\frac{x-y}{x\overline{y}-1}\right|\cdot\left(\frac{x-y}{x\overline{y}-1}\right)^{-1}=\frac{(x\overline{y}-1)|x-y|}{|x\overline{y}-1|(x-y)},
\]
which corresponds to
 \[
a^{-}=\frac{q^{-}y+1}{q^{-}+\overline{y}}=\frac{\left|x-y\right|(x\overline{y}-1) y-(x-y)\left|x\overline{y}-1\right|}{\left|x-y\right|(x\overline{y}-1)-(x-y)\left|x\overline{y}-1\right|\overline{y}}.
\]

\qed

\begin{remark}\label{rem} It follows from this proof that
\begin{eqnarray*}
\frac{x-a^{+}}{y-a^{+}} & = & g(q^{+})=\left(\frac{x\overline{y}-1}{|y|^{2}-1}\right)\cdot\left(\left(\frac{x-y}{x\overline{y}-1}\right) q^{+}+1\right)\\
 & = & \left(\frac{x\overline{y}-1}{|y|^{2}-1}\right)\left(\left(\frac{x-y}{x\overline{y}-1}\right)\frac{(x\overline{y}-1)|x-y|}{|x\overline{y}-1|(x-y)}+1\right)\\
 & = & \left(x\overline{y}-1\right)\left(\frac{|x-y|+|x\overline{y}-1|}{(|y|^{2}-1)|x\overline{y}-1|}\right).
 \end{eqnarray*}
 This observation will be used later.
\end{remark}

\section{Geodesics of weak metrics}\label{s5}

Working in weak metric spaces, it turns out that rather than defining a geodesic as a distance-preserving path (as in the case of metric spaces), it is more
convenient to define it as a path  $\gamma : I \to X$ preserving 
aligned triples (where $I \subset \mathbb{R}$ is some interval). We make the following precise definitions.

\begin{definition}[Aligned triple] Let $(X,\delta)$ be a space with a weak metric and
$x$, $y$ and $z$ be three  points in $X$. We say that the three points 
$x$, $y$, $z$ (in that order) are {\it aligned}
if 
 $\delta(x,z)=\delta(x,y)+\delta(y,z)$.
\end{definition} 
We note that the fact that $x$, $y$, $z$ are aligned does not imply that $z$, $y$, $x$ are aligned.

\begin{definition}[Geodesic] 
A {\it $\delta-$geodesic} (or, simply, a {\it geodesic}) in $X$
is a path $\gamma: I \to X$, where $I$ is an interval of $\mathbb{R}$, 
 such that for any  $t_1$, $t_2$ and $t_3$ in $I$ satisfying $  t_1,\leq t_2\leq t_3 $, the points $\gamma(t_1)$,  $\gamma(t_2)$, $\gamma(t_3)$ are aligned.
\end{definition}

As a simple example, observe that a Euclidean segment is a geodesic
for the Funk weak metric. Observe that in general,  if a path $\gamma:I\to X$ is a  geodesic, then the same path traversed in the opposite direction 
is not necessarily a geodesic.

\medskip

In this section, we discuss geodesics of  Apollonian weak metrics. Let   $A$ be a subset of $\mathbb{E}^n$ satisfying the hypothesis stated at the beginning of
 Section \ref{s4} and let $\delta_A$ be the Apollonian weak metric on $A$. For any $x$ and $y$ in $A$, we consider the following subset of
$\partial A$:
\begin{equation}\label{}
 M_{x,y}= \left\{a_0\in \partial A\ \hbox{ such that } {\vert x-a_0\vert\over \vert y-a_0\vert}= \delta_A(x,y) \right\}.
\end{equation}

We note that in the case where $A$ is bounded, $\partial A$ is compact and nonempty, and therefore, for every $x$ and $y$, $M_{x,y}$ is  nonempty.

\begin{lemma}\label{aligned}
With the above notations, if $x$, $y$ and $z$ are elements in $A$ satisfying
$M_{x,y}\cap M_{y,z} \cap  M_{x,z} \neq  \emptyset$, then the three 
points $x$, $y$ and $z$ are aligned.
\end{lemma}

\proof
This follows from the fact that for any  $a_0$ in $M_{x,z}\cap M_{x,y}\cap M_{y,z}$, we 
have
$$\delta_A(x,y)+\delta_A(y,z)=\log{\vert x-a_0\vert \over \vert y-a_0\vert}+\log{\vert y-a_0\vert\over \vert z-a_0\vert}= 
\log{\vert x-a_0\vert \over \vert z-a_0\vert}=\delta_A(x,z).$$

\qed

Conversely, we have the following

\begin{lemma} Suppose that $A$ is bounded. If $x$, $y$ and $z$ are points in $A$ satisfying
$\delta_A(x,z)=\delta_A(x,y)+\delta_A(y,z)$, then, there exists a point $a_0 \in \partial A$
such that
$$
\delta_A(x,y)={\vert x-a_0\vert\over \vert y-a_0\vert},\
\delta_A(y,z)={\vert y-a_0\vert\over \vert z-a_0\vert} \ \hbox{ and } \
\delta_A(x,z)={\vert x-a_0\vert\over \vert z-a_0\vert}.
$$
\end{lemma}

\proof  Since $\partial A$ is compact and since $z\notin \partial A$, we can find a point $a_0$ in $\partial A$ satisfying $\delta_A(x,z)=\log {\vert x-a_0\vert \over \vert z-a_0\vert}$. This gives
\begin{eqnarray*}
 \log{\vert x-a_0\vert \over \vert z-a_0\vert}  &=&
\log{\vert x-a_0\vert \over \vert y-a_0\vert}+
\log{\vert y-a_0\vert \over \vert z-a_0\vert}  \\ &=&
\sup_{a\in \partial A}\log{\vert x-a\vert \over \vert y-a\vert}+
\sup_{a\in \partial A}\log{\vert y-a\vert \over \vert z-a\vert},
\end{eqnarray*}
which implies  that
$$
\delta_A(x,y) =
\sup_{a\in
\partial A}\log{\vert x-a\vert \over \vert y-a\vert}\ 
=
\log{\vert x-a_0\vert \over \vert y-a_0\vert}
$$
and
$$
\delta_A(y,z) = \sup_{a\in \partial A}\log{\vert y-a\vert \over \vert z-a\vert}
=
\log{\vert y-a_0\vert \over \vert z-a_0\vert}.
$$

\qed

 \medskip

Finally, we apply these results to the case where $A=\mathbb{D}^2$. 
First, we need a lemma about generalized circles in $\mathbb{\overline{C}}=\mathbb{C\cup\mathrm{\left\{ \infty\right\} }}$ (a generalized
circle being as usual a Euclidean circle or a Euclidean straight line compactified by the point $\{\infty\}$). Recall that any generalized circle is the
set images of the unit circle in $\mathbb{C}$ under a M\"obius transformation of $\mathbb{\overline{C}}$.

\medskip

\begin{lemma} \label{lem.critalign}
Four pairwise distinct points $x,y,a,b\in\mathbb{C}$ lie on a generalized
circle if and only if the complex cross-ratio
\[
(x,y,a,b) =\frac{x-a}{y-a}\cdot\frac{y-b}{x-b}
\]
is a real number. \  Furthermore these points appear on the circle in the order $x,y,a,b$ if and only $(x,y,a,b) \in(1,\infty)$.
\end{lemma}

\medskip

\proof This lemma is well known. We recall the proof for the convenience of the reader. 
Let us consider the M\"obius transformation $\phi$ defined by
\[
 \phi(z)= \frac{y-b}{y-a} \cdot \frac{z-a}{z-b}.
\]
The point $x$ belongs to the generalized circle through $a,b,y$
if and only if $\phi(x)$ belongs to the generalized circle through
$\phi(a)=0$, $\phi(b)=\infty$ and $\phi(y)=1$, which is the extended
real line $\mathbb{R}\cup\{\infty\}$. Furthermore, these four points
appear on that circle in the order $x,y,a,b$ if and only if $\phi(x)$
belongs to the interval $(\phi(y),\phi(b))=(1,\infty)$.
The proof of the lemma follows once we observe that $\phi(x) = (x,y,a,b) $.

\qed

 \medskip

Collecting all this information about geodesics of Apollonian weak metrics, we are now ready to
prove the following

\begin{theorem}\label{Th3} 
Let $x$ and $y$ be two distinct points in $\mathbb{D}^{2}$.
Then, the arc of  generalized circle starting at $x$, containing $y$ and orthogonal to the unit circle $\mathbb{S}^{1}$ is a  $\delta_{\mathbb{D}^{2}}$-geodesic starting at $x$ and passing through $y$.
\end{theorem}

\medskip

\proof Denote by $\Gamma$ the generalized circle through $x$ and $y$ and orthogonal to $\mathbb{S}^{1}$. Observe that $\Gamma$ is invariant under the inversion $z \mapsto 1/\overline{z}$,
and therefore $\Gamma$ is  the generalized circle passing through
$x,y,1/\overline{y}$ is orthogonal to the unit circle $\mathbb{S}^{1}$. 

The point $a^{+}(x,y)$ belongs to $\mathbb{S}^{1}$ (see Proposition \ref{S-1}), and Lemma  \ref{lem.critalign} implies that the points $x,y,a^{+}$ appear
in that order on the circle, because $(x,y,a^{+},\frac{1}{\overline{y}}) \in(1,\infty)$.
Indeed, we have from Remark \ref{rem}
\[
\frac{x-a^{+}}{y-a^{+}} = \left(x\overline{y}-1\right)
\left(\frac{|x-y|+|x\overline{y}-1|}{(|y|^{2}-1)|x\overline{y}-1|}\right),
\]
hence
\begin{eqnarray*}
 (x,y,a^{+},\frac{1}{\overline{y}})  & = & \frac{x-a^{+}}{y-a^{+}}\cdot\frac{y-1/\overline{y}}{x-1/\overline{y}}\\
 & = & \left(x\overline{y}-1\right)\left(\frac{|x-y|+|x\overline{y}-1|}{(|y|^{2}-1)|x\overline{y}-1|}\right)\left(\frac{|y|^{2}-1}{x\overline{y}-1}\right)\\
 & = & 1+\frac{|x-y|}{|x\overline{y}-1|}.
 \end{eqnarray*}
We have thus proved that for any pair of points $z,w \in \Gamma$ such that  $z,w,a^{+}$ appear in that order, we have $a^+ \in M_{z,w}$. The Theorem follows now from Lemma \ref{aligned}.

\qed

\end{document}